\documentclass[number,citesort,dvips]{arxbj}
\usepackage{dcolumn}
\usepackage{mathrsfs}
\usepackage{graphicx}

\aid{0}
\volume{17}
\issue{4}
\pubyear{2011}
\firstpage{1400}
\lastpage{1419}
\doi{10.3150/10-BEJ315}

\makeatletter
\newtheorem{lemma}{Lemma}[section]
\newtheorem{theorem}[lemma]{Theorem}
\newtheorem{corollary}[lemma]{Corollary}
\newproclaim{ass}{Assumption}
\newcolumntype{d}[1]{D{.}{.}{#1}}
\makeatother

\begin{document}
\begin{frontmatter}

\title{Estimation for an additive growth curve model with orthogonal
design matrices}
\runtitle{Estimation for an additive growth curve model}
\begin{aug}
\author[a]{\fnms{Jianhua} \snm{Hu}\corref{}\thanksref{a}\ead[label=e1]{frank.jianhuahu@gmail.com}},
\author[b]{\fnms{Guohua} \snm{Yan}\thanksref{b}\ead[label=e2]{gyan@unb.ca}}
\and
\author[a]{\fnms{Jinhong} \snm{You}\thanksref{a}\ead[label=e3]{johnyou07@gmail.com}}
\runauthor{J. Hu, G. Yan and J. You}
\address[a]{School of Statistics and Management, Shanghai
University of Finance and Economics, Shanghai
200433, P.R. China. \printead{e1,e3}}
\address[b]{Department of Mathematics and Statistics, University of New
Brunswick, Fredericton, NB, E3B 5A3 Canada. \printead{e2}}
\end{aug}

\received{\smonth{3} \syear{2008}}
\revised{\smonth{5} \syear{2010}}

%
\begin{abstract}
An additive growth curve model with orthogonal design
matrices is proposed in which observations may have different
profile forms. The proposed model allows us to fit data and then
estimate parameters in a more parsimonious way than the traditional
growth curve model. Two-stage generalized least-squares estimators
for the regression coefficients are derived where a quadratic
estimator for the covariance of observations is taken as the first-stage
estimator. Consistency, asymptotic normality and asymptotic
independence of these estimators are investigated. Simulation
studies and a numerical example are given to illustrate the efficiency
and parsimony of the proposed model for model specifications in the
sense of minimizing Akaike's information criterion (AIC).
\end{abstract}
%
%
\begin{keyword}
\kwd{AIC}
\kwd{asymptotic normality}
\kwd{consistent estimator}
\kwd{growth curve model}
\kwd{quadratic estimator}
\kwd{two-stage generalized least squares}
\end{keyword}

\end{frontmatter}

\section{Introduction}\label{sec1}

In a variety of areas, observations are measured over
multiple time points on a particular characteristic to investigate
the temporal pattern of change on the characteristic. The observations
of repeated measurements are usually analyzed by the growth curve
model (GCM), initiated by Potthoff and Roy \cite{Potthoff64}. Since
then, parameter estimation, hypothesis testing and prediction of
future values have been investigated by numerous researchers,
generating a substantial amount of literature, including \cite{Khatri66,Grizzle69,Rao65,Rao87,Lang88,Lee88,Chaganty03,Kollo07}. The
basic idea of the growth curve model is to
introduce some known functions, usually polynomial functions, so as
to capture patterns of change for time-dependent measurements. We
shall generalize the growth curve model to the case where
observations of time-dependently repeated measurements may have
polynomial functions with different degrees rather than polynomial
functions with a common degree. In this article, different profile
forms mean polynomial functions with different degrees and a profile
form means polynomial functions with a common degree.

To motivate it, let us look at the following situation. We have many
groups of animals, with each group being subjected to a different
treatment. Animals in all groups are measured at the same $p$
time points and assumed to have the same covariance matrix
$\Sigma$. The growth curve associated with the $i$th group is $\theta
_{i0}+\theta_{i1}t+\theta_{i2}t^2+\cdots+\theta_{iq_i}t^{q_i}$,
implying that the growth curves may have different profiles, say $k$
profiles, not necessarily one profile. There are $m_i$ groups that
have the same profile form with index $i$ and $n_i$
individuals in total. Here $n=\sum_{i=1}^kn_i$. The
simplest situation is that each group has a different profile form.
Assume that there are $k$ groups of individuals and $p$ observing
time points such that $k+p\leq n$. For $i=1,2,\ldots,k$, put
\begin{eqnarray*}
{Z}_i&=&\left[\matrix{
1 & t_1 & t_1^2 & \ldots& t_1^{q_i-1} \vspace*{2pt}\cr
1 & t_2 & t_2^2 & \ldots& t_2^{q_i-1} \vspace*{2pt}\cr
. & . & . & \ldots& . \vspace*{2pt}\cr
1 & t_p & t_p^2 & \ldots& t_p^{q_i-1}}\right] ,
\\
{\Theta}_i&=&(\theta_{i0},\theta_{i1},\theta_{i2},\ldots,\theta
_{i q_i-1}),
\end{eqnarray*}
and
\[
{X}_i=(x_{i1},x_{i2},\ldots,x_{in})^{\prime}\in\mathbb{R}^n,
\]
where $x_{i(p_i+j)}=1$ for $j=1,2,\ldots,n_i$ with $p_0=0$, $%
p_i=\sum_{j=1}^{i-1}n_j$ and other $x^{\prime}_{ij}$'s are $0$.

Generalizing the above situation we propose the following additive
growth curve model
\begin{equation}\label{e11}
Y=\sum_{i=1}^kX_i\Theta_i Z_i' +\mathcal{E},\qquad \mathcal{E}\sim
\mathcal{G}(\mathbf{0},I \otimes\Sigma),
\end{equation}
with orthogonal design matrices or mutually orthogonal column spaces
of design matrices, defined as
\begin{equation}\label{e12}
\operatorname{rank}(X_i)+p\leq n\quad  \mbox{and}\quad  \mathscr C(X_i)\perp
\mathscr
C(X_j) \quad \mbox{or}\quad  X_i'X_j=\mathbf{0} \mbox{ for any distinct } i,j,
\end{equation}
where $Y $ is an $n \times p$ matrix of observations; $X_i$, $Z_i$
$(1\leq i \leq k)$ are known $n\times m_i$ ($n>m_i$) full-rank
design matrices and $p\times q_i$ ($p>q_i$) full-rank profile
matrices, respectively; $\Theta_i$ ($1\leq i \leq k$) are unknown
$m_i\times q_i$ matrices of the regression coefficients; $\mathscr
C(X)$ denotes the column space of the matrix $X$; $\mathcal{G}$ is a
general continuous type distribution function; observations on
individuals are independent; and the rows of the random error matrix
$\mathcal E$ are independent and identically distributed with mean
zero and a common unknown covariance matrix $\Sigma$ of order $p$.

The model (\ref{e11}) subject to (\ref{e12}) will be demonstrated
to have an advantage that it fits data in a more parsimonious way than
the traditional growth curve model in the situation
where model specification is needed. In the above stated example of
animals, the traditional growth curve model assumes that all
observations have the same profile form, which may cause the
model misspecification, underfitting or overfitting.

On the other hand, Kollo and von Rosen \cite{Kollo05}, in Chapter 4,
investigated an additive growth curve model with nested column
spaces generalized by design matrices, that is, constraint $\mathscr
C(X_1)\supseteq\mathscr C(X_2)\supseteq\cdots\supseteq\mathscr C(X_k)$
with $\operatorname{rank}(X_1)+p\leq n$, usually called the extended growth
curve model. Obviously, there is not an inclusion relationship between
the extended growth curve model and the proposed model (\ref{e11})
with (\ref{e12}) because the constraint of nested column spaces and
the constraint of orthogonal column spaces have no inclusion
relationship. An extension of the growth curve model proposed in
\cite{Verbyla88} did not include the proposed
model~(\ref{e11}) with (\ref{e12}), either.

This paper will investigate estimation of parameters and
properties of the corresponding estimators in the proposed model
(\ref{e11}) with (\ref{e12}), including consistency and asymptotic
normality.


The organization of the paper is as follows. Two-stage generalized
least-squares estimators of the regression coefficients are obtained
in Section \ref{sec2}. Both the consistency of the estimators for the regression
coefficients and a quadratic estimator for the unknown covariance are
investigated in Section \ref{sec3}, while their asymptotic normalities under
certain conditions are
investigated in Section~\ref{sec4}. Simulation
studies are given in Section \ref{sec5}. A numerical example is explored to
illustrate our techniques in Section \ref{sec6}. Finally, brief
concluding remarks are stated in Section \ref{sec7}.

Throughout this paper, the following notations are used. $\mathscr
M_{n\times p}$ denotes the set of all $n\times p$ matrices over real
set $\mathbb{R}$ with trace inner product $\langle,\rangle$ and $\|\cdot\|$ denotes
the trace norm on the set $\mathscr M_{n\times p}$. $\operatorname{tr}(A)$
denotes the trace of matrix $A $ and $I_n $ denotes the identity
matrix of order $n$. For an $n\times p$ matrix $Y$, we write $Y
=[\mathbf{y}_1', \ldots, \mathbf{y}_n']'$, $\mathbf{y}_i'\in
\mathbb
{R}^p$, where $\mathbb{R}^p$ is the $p$-dimensional
real space and $\operatorname{vec}(Y)$ denotes $np$-dimensional vector
$[\mathbf{y}_1, \ldots, \mathbf{y}_n]'$. Here the $\operatorname
{vec}$ operator
transforms a matrix into a vector by stacking the rows of the matrix
one underneath another. $Y\sim\mathcal{G}(\bolds{\mu},I\otimes
\Sigma)$
means that $Y$ follows a general continuous type distribution
$\mathcal{G}$ with $\mathrm{E}(Y)=\bolds{\mu}$ and $ \mathrm
{E}(Y-\bolds{\mu})(Y-\bolds{\mu})'=I\otimes
\Sigma$. The Kronecker product $A \otimes B $ of matrices $A $
and $B $ is defined to be $A \otimes B =(a_{ij}B )$. Then we
have $\operatorname{vec}(A B C)=(A \otimes C')\operatorname{vec}(B
)$. Let $A ^+$ denote
the Moore--Penrose inverse of $A $. $P_X =X (X 'X )^-X '$ denotes the
orthogonal projection onto the column space $\mathscr C(X) $ of a
matrix~$X $. $M_{X}=I-X(X'X)^-X'$ is the orthogonal projection onto
the orthogonal complement $\mathscr C(X)^\perp$ of $\mathscr C(X)$.

\section{Two-stage generalized least-squares estimators}\label{sec2}

Recall that the regression coefficients, $\Theta_1,\ldots,\Theta_k$,
in the model (\ref{e11}) are defined before a design is planned and
observation $Y$ is obtained. Thus the rows of the design matrices,
$X_1,\ldots,X_k$, are added one after another and the profile forms,
$Z_1,\ldots,Z_k$, do not depend on the sample size $n$. So, we shall
only consider the case of full-rank $X_i$'s and $Z_i $'s in the
present paper.

Set
\begin{equation}\label{e21}
\bolds{\mu}=\sum_{i=1}^kX_i\Theta_i Z_i'.
\end{equation}
Equation (\ref{e21}) is said to be the mean structure of the model
$(\ref
{e11})$.

A statistic $\hat{\bolds{\mu}}_{\mathrm{gls}}(Y)$ is said to be the generalized
least-squares (GLS) estimator of parameter matrix $\bolds{\mu}$ if the
minimum value of function $\langle Y-\bolds{\mu},Y-\bolds{\mu}\rangle$ is
attained at the point
$\mathbf{u}=\hat{\bolds{\mu}}_{\mathrm{gls}}(Y)$, where the inner product
$\langle,\rangle$ or
the trace norm $\|\cdot\|$ associated with the covariance $I\otimes\Sigma$
of $Y\dvtx \langle\mathbf{w}_1,\mathbf{w}_2\rangle=\operatorname{vec}(\mathbf{w}_2)'(I
\otimes
\Sigma)^{-1}\operatorname{vec}(\mathbf{w}_1)$ with $\parallel
\mathbf{w} \parallel
=\langle\mathbf{w},\mathbf{w}\rangle^{1/2}$ and $\mathbf{w},\mathbf
{w}_1,\mathbf
{w}_2\in\mathscr M_{n\times p}$.

Generally speaking, we actually know nothing or very little about
the covariance $\Sigma$ of observations of repeated measurements
before we measure these observations. So, alternatively, a~two-stage
estimation is used to find an estimator of $\bolds{\mu}$, denoted by
$\hat{\bolds{\mu}}_{\mathrm{2sgls}}(Y)$. The two-stage estimation procedure
is as
follows: First, based on the observation $Y $, find a first-stage
estimator $\widetilde{\Sigma}$ of $\Sigma$. Second, replace
the unknown $\Sigma$ with the first-stage estimator
$\widetilde{\Sigma}$ and then find $\widehat{\bolds{\mu}}_{\mathrm{2sgls}}(Y)$ through
the GLS method. For convenience, we shall omit
the subscript of $\hat{\bolds{\mu}}_{\mathrm{2sgls}}(Y)$.

In order to get a good first-stage estimator $\widetilde{\Sigma}$ of
$\Sigma$, let us have a close look at the following quadratic
statistic (a quadratic form without associating with parameters):
\begin{equation}\label{e22}
\widehat{\Sigma}(Y )=Y'W Y, \qquad W \equiv\frac
{1}{r}\Biggl(I-\sum_{i=1}^kP_{X_i }\Biggr),
\end{equation}
where $r=n-\sum_{i=1}^k\operatorname{rank}(X_i)$.
\begin{longlist}[(3)]
\item[(1)] The statistic\vspace*{1pt} $\widehat{\Sigma}(Y)$ is easily proven to be
positive definite with probability $1$; see Theorem~3.1.4 of
\cite{Muirhead82}. So, $\widehat{\Sigma}^{-1}(Y)$ exists with
probability 1.

\item[(2)] Under the assumption of normality, the quadratic estimator
$\widehat{\Sigma}(Y)$ given by equation (\ref{e22}) follows a Wishart
distribution;
see \cite{Hu08}.

\item[(3)] $\widehat\Sigma(Y )$ is an unbiased invariant estimator of $\Sigma$; see \cite{Hu10}. A similar result for the
growth curve model was obtained by \v{Z}e\v{z}ula \cite{Zezula93}.
\end{longlist}

It follows from the above properties that the statistic
$\widehat\Sigma(Y)$ seems to be a very good candidate for the
first-stage estimator. As a consequence, it will be taken as the
first-stage estimator $\widetilde{\Sigma}$ in our subsequent
discussion.

For $ i=1,\ldots,k$, let
\begin{equation}\label{e25}
H_i (Y)\equiv\widehat\Sigma^{-1}(Y )Z_i(Z_i'
\widehat\Sigma^{-1}(Y)Z_i )^{-1}Z_i'=\widehat\Sigma^{-1}(Y
)(P_{Z_i} \widehat\Sigma^{-1}(Y)P_{Z_i} )^+.
\end{equation}
Then, we easily see that
\begin{equation}\label{e26}
Z_i'H_i (Y)=Z_i', \qquad i=1,\ldots,k.
\end{equation}

When $\widehat\Sigma(Y)$ is taken as the first-stage estimator, the
following lemma provides the explicit expression of the two-stage
GLS estimators both for mean matrix $\bolds{\mu}$ and
the regression coefficients, $\Theta_1, \ldots, \Theta_k$.
Furthermore, under certain conditions, these estimators are
unbiased.

\begin{theorem}\label{l21} Consider $\Sigma=\widehat\Sigma(Y)$ for the
model (\ref{e11}) subject to (\ref
{e12}). The following statements hold.
\begin{enumerate}[(1)]
\item[(1)] The two-stage GLS estimator
$\hat{\bolds{\mu}}(Y)$ of $\bolds{\mu}$ is given by
\begin{equation}\label{e27}
\hat{\bolds{\mu}}(Y)=\sum_{i=1}^kP_{X_i}Y\Sigma
^{-1}(Y)(P_{Z_i}\Sigma
^{-1}(Y)P_{Z_i})^+=\sum_{i=1}^k
P_{X_i}YH_i (Y ).
\end{equation}
\item[(2)] The two-stage GLS estimator
$\widehat\Theta_i(Y)$ of $\Theta_i$ is given by
\begin{equation}\label{e28}
\widehat\Theta_i(Y)=(X_i'X_i)^{-1}X_i'YH_i (Y)Z_i(Z_i'Z_i)^{-1}.
\end{equation}
\item[(3)] If the distribution of $\mathcal E $ is symmetric about the origin
$\mathbf{0}$, the statistic $\hat{\bolds{\mu}}(Y)$ is an unbiased
estimator of mean $\bolds{\mu}$. Moreover, for each $i$, the statistic
$\widehat\Theta_i(Y)$ is an unbiased estimator of the regression
coefficients $\Theta_i$.
\end{enumerate}
\end{theorem}

The proof of Theorem \ref{l21} is deferred to the \hyperref[app]{Appendix}.

\section{Consistency}\label{sec3}

Since $Y $ is associated with the sample size $n$, we shall use $Y
_n$ to replace $Y$ in (\ref{e22})--(\ref{e28}) and then investigate
the consistency of the estimator $\widehat\Sigma(Y _n)$ and the
estimators, $\widehat\Theta_1(Y),\ldots,\widehat\Theta_k(Y)$, as the
sample size $n$ tends to infinity. Note that $X$ and $\mathcal E$
are also associated with the sample size $n$.

Regarding the consistency of the quadratic estimator
$\widehat\Sigma(Y_n)$, we have the following result.
\begin{theorem}\label{t31}
For the model (\ref{e11}) subject to (\ref{e12}), the statistic
$\widehat\Sigma(Y_n)$ defined by (\ref{e22}) is a consistent
estimator of the covariance matrix $\Sigma$.
\end{theorem}

\begin{pf} It follows from the invariance of statistic
$\widehat\Sigma(Y)$ that $\widehat\Sigma(Y)=\widehat\Sigma
(\mathcal
E)$. And $\widehat\Sigma(Y)$ can be rewritten as
\begin{equation}\label{e31}
\widehat{\Sigma}(Y_n)=\widehat{\Sigma}(\mathcal
E)=\frac{n}{n-m}\Biggl(\frac{1}{n}\sum_{l=1}^n\mathcal{E}_l\mathcal{E}_l'-\frac{1}{n}\mathcal{E}'\sum_{i=1}^kP_{X_i}\mathcal{E}\Biggr),
\end{equation}
where $m=\sum_{i=1}^k\operatorname{rank}(X_i)$ and $\mathcal
E=(\mathcal
E_1,\ldots,\mathcal E_n)'\sim\mathcal{G}(\mathbf{0}, I_n \otimes
\Sigma)$.

Note that $(\mathcal E_l\mathcal E_l')_{l=1}^n$ is a random sample
from a population with mean $\mathrm{E}(\mathcal E_l \mathcal
E_l')=\Sigma$.
Kolmogorov's strong law of large numbers tells us that
\begin{equation}\label{e32}
\frac{1}{n}\sum_{l=1}^n\mathcal E_l\mathcal E_l'
\mbox{ converges almost surely to } \Sigma.
\end{equation}

Let $\varepsilon>0$. By Chebyshev's inequality and $\mathrm
{E}(\mathcal
E'\mathcal E)=\operatorname{tr}(I )\Sigma$, we have
\begin{eqnarray*}
P\Biggl(\Biggl\Vert\frac{1}{\sqrt{n}}\sum_{i=1}^kP_{X_i}\mathcal{E}\Biggr\Vert\geq\varepsilon\Biggr)
& \leq&\frac{1}{n\varepsilon
^2}\mathrm{E}\Biggl[\operatorname{tr}\Biggl(\mathcal{E}'\sum_{i=1}^kP_{X_i}\mathcal{E}\Biggr)\Biggr] =
\frac{1}{n\varepsilon^2}\operatorname{tr}\Biggl(\mathrm{E}[\mathcal{E}\mathcal{E^{'}}]\sum_{i=1}^kP_{X_i}\Biggr)\\
&=&\frac{1}{n\varepsilon^2}\operatorname{tr}\Biggl(I_n\operatorname
{tr}(\Sigma)\sum_{i=1}^k
P_{X_i}\Biggr)=\frac{1}{n\varepsilon
^2}\operatorname{tr}\Biggl(\sum_{i=1}^kP_{X_i}\Biggr)\operatorname{tr}(\Sigma).
\end{eqnarray*}
Since $ \operatorname{tr}(\sum_{i=1}^kP_{X_i})=\sum
_{i=1}^k\operatorname{rank}(X_i)$ is a
constant, $P(\parallel\frac{1}{\sqrt{n}}\sum_{i=1}^kP_{X_i}\mathcal{E}\parallel\geq\varepsilon)$ tends to $\mathbf{0}$ as the sample
size $n$ tends to infinity. So
\begin{equation}\label{e33}
\frac{1}{\sqrt{n}}\sum_{i=1}^kP_{X_i}\mathcal{E} \mbox{ converges
in probability to } \mathbf{0}.
\end{equation}
Since convergence almost surely implies convergence in probability,
by (\ref{e32}) and (\ref{e33}), we obtain from (\ref{e31}) that
$\hat{\Sigma}(Y _n)$ converges in probability to $\Sigma$, which
completes the proof.
\end{pf}

\begin{ass}\label{ass1} For $l=1,\ldots,k$,
\begin{equation}\label{e34}
\lim_{n\rightarrow\infty}n^{-1}X _l'X_l =R_l,
\end{equation}
where $R_l$ is positive definite.
\end{ass}

For convenience, we restate Lemma 3.2 of \cite{HuYan08} as follows.
\begin{lemma}\label{l31}
For $i\in\{1,\ldots,k\}$, $H_i(Y_n)$ converges in probability to
$H_i\equiv\Sigma^{-1}(P_{Z_i}\Sigma^{-1}P_{Z_i})^+.$
\end{lemma}

On the consistency of the estimators of the regression coefficients
$\Theta_i(Y _n)$s, we obtain the following theorem.
\begin{theorem}\label{t32}
For any fixed $i\in\{1,\ldots,k\}$, with Assumption \ref{ass1}, the statistic
$\widehat\Theta_i(Y_n)$ is a consistent estimator of the
regression coefficient $\Theta_i$.
\end{theorem}

\begin{pf} Fix $i\in\{1,\ldots,k\}$. By equation (\ref{e28}),
we obtain the following equation:
\begin{equation}\label{e35}
\widehat\Theta_i(Y_n)=\Theta_i+S_i\mathcal E H_i(Y_n)K_i,
\end{equation}
where $S_i=(X_i'X_i)^{-1}X_i'$ and $K_i=Z_i(Z_i'Z_i)^{-1}$. The
second term of the right side in (\ref{e35}) can be rewritten as
\[
S_i\mathcal EH_i(Y_n)K_i=n(X_i'X_i)^{-1}
\biggl(\frac{1}{\sqrt{n}}X_i'\biggr)\biggl(\frac{1}{\sqrt{n}}P_{X_i}
\mathcal{E}\biggr)H_i(Y _n)K_i.
\]
By condition (\ref{e34}), $X_i'/{\sqrt{n}}$ are bounded. In fact,
the elements of $X_i'/\sqrt{n}$ are at most of order $n^{-1/2}$ (see
the proof of Lemma {\ref{l41}} below). So by (\ref{e33}), (\ref
{e34}), Lemma \ref{l31} and Theorem 11.2.12 of \cite{Lehmann05},
the second term of the right side in (\ref{e35})
converges in probability to~$\mathbf{0}$. Thus, $\widehat\Theta
_i(Y_n)$ converges in probability to $\Theta_i$, which completes the
proof.
\end{pf}

In order to prove the consistency of the estimators $\widehat\Theta
_1(Y _n)$, $\widehat\Theta_2(Y _n), \ldots, \widehat\Theta_k(Y
_n)$, the conditions $\lim_{n\rightarrow\infty}n^{-1}X _l'X_l =R_l$
for $l=1,\ldots,k$ have been used in Theorem \ref{t32}. We imagine
that for each new observation, a new row is added to the matrices
$X_l$ and that the earlier rows remain intact in such a way that,
for $l=1,2,\ldots,k$, the $m_l\times m_l$ elements of $X _l'X_l $
are $\mathrm{O}(n)$. In addition, we exclude the possibility that the limits
of $n^{-1}X _l'X_l$s are singular.

\section{Asymptotic normality}\label{sec4}

We have investigated the consistency of the estimators
$\widehat\Sigma(Y)$ and $\widehat\Theta_i(Y_n)$ in the preceding
section. In this section, we shall investigate the asymptotic
normality of $\sqrt{n}[\widehat\Theta_i(Y_n)-\Theta_i]$ and
$\sqrt{n}[\widehat\Sigma(Y_n)-\Sigma]$ under certain
conditions.

We need the following lemma in the proof of the subsequent results.

\begin{lemma}\label{l41} Let $S_i$ $=$ $(X_i 'X_i )^{-1}X_i '$ $\equiv
(\mathbf{s}_{i1},$ $\mathbf{s}_{i2}, \ldots,\mathbf{s}
_{in})_{m_i\times n}$,
where $\mathbf{s}_{ij}$ is the $j$th column of $X_i$. Then, under
condition (\ref{e34}), the $m_i$ elements of $\sqrt{n}\mathbf{s}_{ij}$
are $\mathrm{O}(n^{-1/2})$ for any $i\in\{1,\ldots,k\}$ and
$j\in\{1,\ldots,n\}$.
\end{lemma}

The proof of Lemma \ref{l41} is deferred to the \hyperref[app]{Appendix}.

\begin{theorem}\label{t41} Under Assumption \ref{ass1}, the random matrix
$\sqrt
{n}S_i\mathcal E$ converges in
distribution to $\mathcal N_{m_i\times p}(\mathbf{0}, R_i^{-1}\otimes
\Sigma)$ for any $i\in\{1,\ldots,k\}$.
\end{theorem}

Also, the proof of Theorem \ref{t41} is deferred to the \hyperref[app]{Appendix}.

Finally, by Theorem \ref{t41} and Slutsky's theorem, we obtain our
main result on the asymptotic normality of
$\sqrt{n}[\widehat\Theta_i(Y_n)-\Theta_i]$.

\begin{theorem}\label{t42} Under Assumption \ref{ass1}, the statistic $\sqrt
{n}[\widehat\Theta_i(Y _n)-\Theta_i]$
converges in distribution to $\mathcal N_{m_i\times q_i}(\mathbf
{0},R_i\otimes(Z_i'\Sigma Z_i)^{-1})$ for any $i\in\{1,\ldots,k\}$.
\end{theorem}

Next, we shall investigate the asymptotic normality of the
$\widehat\Sigma(Y)$. The fourth-order moment of the error matrix
will be needed in the following discussion.

\begin{ass}\label{ass2} $\mathrm{E}(\mathcal E_1)=\mathbf{0}$,
$\mathrm{E}(\mathcal E_1\mathcal E_1')=\Sigma>\mathbf{0}$,
$\mathrm{E}(\mathcal E_1\otimes\mathcal E_1\mathcal E_1')=\mathbf
{0}_{p^2\times p}$ and $\mathrm{E}\Vert\mathcal E_1\Vert ^4<\infty$, where
$\mathcal E_1'$ is the first row vector of the error matrix
$\mathcal E$.
\end{ass}

\begin{theorem}\label{t43}
Under Assumptions \ref{ass1} and \ref{ass2}, the following probability
statements hold:
\begin{longlist}
\item[(a)] $\sqrt{n}(\widehat\Sigma(Y)-\Sigma)$ converges to
$\mathcal N(\mathbf{0}, \operatorname{Cov}(\mathcal E_1'\otimes
\mathcal
E_1'))$ in distribution.

\item[(b)] For each $i$, $\sqrt{n}(\widehat\Sigma(Y)-\Sigma)$
and $\sqrt{n}(\widehat\Theta_i(Y)-\Theta_i)$ are
asymptotically independent.

\item[(c)] For any distinct $i$, $j$, $\sqrt{n}(\widehat\Theta
_i(Y)-\Theta_i)$ and $\sqrt{n}(\widehat\Theta
_j(Y)-\Theta_j)$ are independent.
\end{longlist}
\end{theorem}

\begin{pf} (a)
$\sqrt{n}(\widehat\Sigma(Y)-\Sigma)$ can be decomposed
into
\[
\sqrt{n}\bigl(\widehat\Sigma(Y)-\Sigma\bigr)=\Delta_1+\Delta_2+\Delta_3,
\]
where
\begin{eqnarray*}
\Delta_1&=&\sqrt{n}\Biggl(\frac{1}{n}\sum_{l=1}^n\mathcal E_l\mathcal
E_l'-\Sigma\Biggr),
\\
\Delta_2&=&\frac{m}{\sqrt{n}(n-m)}\sum_{l=1}^n\mathcal E_l\mathcal
E_l',
\\
\Delta_3&=&-\frac{\sqrt{n}}{n-m}\mathcal
E'\sum_{i=1}^kP_{X_i}\mathcal E.
\end{eqnarray*}

Similar to the proof of conclusions (\ref{e32}) and (\ref{e33}) in
Theorem \ref{t31}, we easily obtain that $\Delta_2$ and~$\Delta_3$
converges to $\mathbf{0}$ in probability 1.

Also by assumptions 1 and 2, $\Delta_1$ converges to
$\mathcal N(\mathbf{0}, \Phi_2)$ in distribution, where
$\Phi_2=\operatorname{Cov}(\mathcal E_1'\otimes\mathcal E_1')$.
Thus, we
have
\[
\sqrt{n}\operatorname{vec}\bigl(\widehat\Sigma(Y)-\Sigma\bigr)=\operatorname
{vec}(\Delta_1)+\mathrm{o}_P(\mathbf{1}).
\]
Hence, $\sqrt{n}(\widehat\Sigma(Y)-\Sigma)$ converges to
$\mathcal N(\mathbf{0}, \operatorname{Cov}(\mathcal E_1'\otimes
\mathcal
E_1'))$ in distribution.

(b) By equation (\ref{e35}), it suffices to prove the asymptotic
independence between $\frac1{\sqrt{n}}\operatorname
{vec}(X_i'\mathcal E)$
and $\sqrt{n}\operatorname{vec}(\widehat\Sigma(Y)-\Sigma)$.

Let $Q_n=X_i'\mathcal E=(\mathbf{x}_1^i,\ldots,\mathbf
{x}_n^i)(\mathcal
E_1,\ldots,\mathcal E_n)'$. Then
\begin{eqnarray*}
\operatorname{Cov}\biggl(\biggl(\frac1{\sqrt{n}}X_i'\mathcal
E\biggr),\sqrt{n}(\widehat\Sigma-\Sigma)\biggr) &=&
\operatorname{Cov}\Biggl(\Biggl( \sum_{l=1}^n\mathbf{x}_l^i\mathcal
E_l'\Biggr),\Biggl(\frac1n \sum_{l=1}^n\mathcal
E_l\mathcal E_l'
-\Sigma\Biggr)\Biggr)+\mathrm{o}_P(\mathbf1)\\
&=&\mathrm{E}\Biggl(\Biggl( \sum_{l=1}^n\mathbf{x}_l^i\otimes\mathcal
E_l'\Biggr)\Biggl(\frac
1n \sum_{j=1}^n\mathcal E_j\otimes\mathcal
E_j'-\Sigma\Biggr)\Biggr)+\mathrm{o}_P(\mathbf1).
\end{eqnarray*}
According to Assumption \ref{ass2}, $\operatorname{Cov}((\frac
1{\sqrt{n}}X_i'\mathcal
E),\sqrt{n}(\widehat\Sigma(Y)-\Sigma))$
converges to $\mathbf{0}$ in probability 1. It follows that the
vectors $\frac1{\sqrt{n}}\operatorname{vec}(X_i'\mathcal E)$
and $\sqrt{n}\operatorname{vec}(\widehat\Sigma(Y)-\Sigma)$ are
asymptotically independent. Therefore,
$\sqrt{n}(\widehat\Sigma(Y)-\Sigma)$ and
$\sqrt{n}(\widehat\Theta_i(Y)-\Theta_i)$ also are
asymptotically independent.

(c) For any distinct $i,j$, it follows from condition (\ref
{e12}) that
\[
\operatorname{Cov}\bigl(\sqrt{n}\bigl(\widehat\Theta_i(Y)-\Theta_i\bigr),
\sqrt{n}\bigl(\widehat\Theta_j(Y)-\Theta_j\bigr) \bigr)=\mathbf{0}.
\]

We have completed the proofs of statements (A)--(C).
\end{pf}

Sometimes, it is necessary to consider hypothesis tests of the
form
\[
\mbox{H}_i\dvtx C\Theta_iV'=\mathbf{0},
\]
where $C$ and $V$ are, respectively, $s\times m_i$ and $t\times q_i$
constant matrices. In this case, Theorem \ref{t42} and Slutsky's
theorem are explored to understand the asymptotic behavior of
$\sqrt{n}(C\widehat\Theta_i(Y)V'-C\Theta_iV')$.

\begin{corollary}\label{c41}
Under Assumption \ref{ass1}, if matrices $C (X'X)^{-1}C'$ and $V(
Z'\widehat\Sigma^{-1}(Y)Z)^{-1}V'$ are non-singular, then the
statistic
\[
(Cn(X'X)^{-1}C')^{-1/2}\sqrt{n}(C\widehat\Theta_i(Y)
V')(V(Z'\widehat\Sigma^{-1}(Y)Z)^{-1}
V')^{-1/2}
\]
under H$_i$ converges in distribution to $\mathcal N_{s\times
t}(\mathbf{0},I)$.
\end{corollary}

Therefore, if it is necessary to make the statistical inference about
certain $\Theta_i$ in the model (\ref{e11}), we can take the
normal distribution $\mathcal N_{s\times t}(0,I)$ as an approximate
distribution of the statistic
\[
(Cn(
X'X)^{-1}C')^{-1/2}\sqrt{n}(C\widehat\Theta
_i(Y)V')(V(Z'\widehat
\Sigma^{-1}(Y)Z)^{-1}V')^{-1/2}
\]
under H$_i$ if the sample size is large. Moreover, due to (c) of
Theorem \ref{t43},
H$_i$ and H$_j$ can be considered independently.

\section{Simulation studies}\label{sec5}

In this section, we shall use simulation to investigate the
efficiency and parsimony of the model~(\ref{e11}) subject to
constraint (\ref{e12}), compared with the traditional growth curve
model $Y=X\Theta Z'+\mathcal E$.

We take an example as follows. Suppose $n$ patients are divided into two
groups with numbers of patients $n_1$ and $n_2$, respectively. A certain
measurement in an active drug trial is made on each of the $n_1$ patients
taking a placebo and the $n_2$ patients taking the active drug at time
points $t_1=-1$, $t_2=-0.5$, $t_3=0.5$ and $t_4=1$. Assume that the
first $n_1$ observations come from normal distribution $\mathcal
N(\bolds{\mu}_1, \Sigma_0),$ where
\[
\bolds{\mu}_1=(4+2t_1,4+2t_2, 4+2t_3, 4+2t_4)
\]
and
\[
\Sigma_0=
\pmatrix{
1&\rho&\rho^2 &\rho^3\cr
\rho& 1 & \rho&
\rho^2\cr
\rho^2 & \rho& 1 & \rho\cr
\rho^3 & \rho^2 & \rho&
1}
,\qquad 0< \rho< 1.
\]
(The model with this $\Sigma_0$ is called the simplest serial
correlation model in literature.) It means that the $n_1$
observations have a linear profile form of time points. The
remaining $n_2$ observations come from normal distribution $\mathcal
N(\bolds{\mu}_2, \Sigma_0),$ where
\[
{\bolds{\mu}_2}=(3+2t_1+t_1^2-t_1^3, 3+2t_2+t_2^2-t_2^3,
3+2t_3+t_3^2-t_3^3, 3+2t_4+t_4^2-t_4^3),
\]
implying that the $n_2$ observations have a cubic polynomial
profile form of time points.

Let
\[
Z_1'=
\pmatrix{
1&1&1& 1\cr
t_1 & t_2 &t_3& t_4}
\quad \mbox{and}\quad Z_2'=
\pmatrix{
1&1&1& 1\cr
t_1 & t_2 &t_3& t_4 \vspace*{2pt}\cr
t_1^{2}& t_2^{2}&t_3^2& t_4^{2}\vspace*{2pt}\cr
t_1^{3}& t_2^{3}&t_3^3&
t_4^{3}}
\]
and
\[
B_1=
\pmatrix{ 4 & 2},\qquad B_2=
\pmatrix{ 3 & 2& -3
& 2
}
,
\]
then
\[
\bolds{\mu}_1=B_1Z_1'\quad \mbox{and}\quad  \bolds{\mu}_2=B_2Z_2'.
\]

In real experiments, with observations $Y$, model specification is a
challenging task.

Without loss of generality, we shall consider three approaches using
the growth curve model to fit data of repeated
measurements from the above synthetic example.

The first approach is to regard all observations of repeated
measurements as having linear profile forms over multiple time points.
In this scenario, model underfitting has occurred. The underfitted
model is denoted by
$\psi_{\mathrm{u}}$,
\[
\operatorname{Model} \psi_{\mathrm{u}}\dvtx Y=X\Theta_uZ_1'+\mathcal E,
\]
where\vspace*{1pt} $X=
\bigl({\mathbf1_{n_1}\enskip\mathbf{0}\atop
\mathbf{0} \enskip\mathbf1_{n_2}}\bigr)
$ and
$\Theta_u=
\bigl({\theta_{11}^u \enskip \theta_{12}^u\atop
\theta_{21}^u\enskip\theta_{22}^u}\bigr)
$ to fit the $n$
observations. By (b) of Lemma \ref{l21}, the estimator
$\widehat\Theta
_u=(X'X)^{-1}X'Y\widehat\Sigma^{-1}(Y)Z_1(Z_1'\widehat\Sigma
^{-1}(Y)Z_1)^{-1}$.

The second approach is to regard all observations of repeated
measurements as following cubic polynomial profile forms over multiple
time points. In this case, model overfitting has occurred. The
overfitted model is denoted by
$\psi_{\mathrm{o}}$,
\[
\operatorname{Model} \psi_{\mathrm{o}}\dvtx Y=X\Theta_oZ_2'+\mathcal E,
\]
where\vspace*{1pt} $\Theta_o=
\bigl({\theta_{11}^o \enskip \theta_{12}^o\enskip
\theta_{13}^o \enskip \theta_{14}^o\atop \theta_{21}^o\enskip
\theta_{22}^o\enskip
\theta_{23}^o\enskip\theta_{24}^o
}\bigr)$ to fit the $n$
observations. The estimators of regression coefficients are
$\widehat\Theta
_o=(X'X)^{-1}X'Y\widehat\Sigma^{-1}(Y)Z_2(Z_2'\widehat\Sigma
^{-1}(Y)Z_2)^{-1}$.

The third approach is to regard the first $n_1$ observations as having
a linear profile form of time points and the remaining $n_2$
observations as having a cubic polynomial profile form over multiple
time points. In this case, model misspecification may not occur.
The additive model is denoted by
$\psi_{\mathrm{a}}$,
\[
\operatorname{Model} \psi_{\mathrm{a}}\dvtx Y=X_1\Theta
_1Z_1'+X_2\Theta
_2Z_2'+\mathcal E,
\]
where $X_1=
\bigl({\mathbf1_{n_1}\atop \mathbf{0}}
\bigr)
$,
$X_2=
\bigl({\mathbf{0}\atop \mathbf1_{n_2}}
\bigr)
$,
$\Theta_1=
\bigl({\theta_{11} \enskip \theta_{12}}
\bigr)
$
and $\Theta_2=
\bigl({\theta_{21} \enskip\theta_{22}\enskip
\theta_{23}\enskip \theta_{24}}
\bigr)$, to fit the $n$ observations.
Based on the above assumption, $\psi_{\mathrm{a}}$ actually is the
true model. The estimators $\widehat\Theta
_i=(X_i'X_i)^{-1}X_i'Y\widehat\Sigma^{-1}(Y)Z_i(Z_i'\widehat\Sigma
^{-1}(Y)Z_i)^{-1}$
for $i=1,2$.

The model specification starts with residuals. We shall use a residual matrix
$R$, defined as the difference between the observation $Y$ and
fitted mean $\widehat Y $, that is, $R=Y-\widehat Y$, to discuss the
model specification for our example.

The residual matrix sum of squares (RMSS) is defined as the trace of
$R'R$
\begin{equation}\label{e51}
\operatorname{RMSS}=\Vert R\Vert ^2=\operatorname{tr}\bigl((Y-\widehat Y
)'(Y-\widehat Y )\bigr).
\end{equation}

Usually, overfitting of model specification can provide a smaller RMSS as
well as use more parameters. The residual matrix sum of squares and the
number of
parameters are two trade-off issues in model specification. Here,
Akaike's information criterion (AIC) -- see \cite{Akaike73} -- is explored
to reward the decreasing RMSS and penalize overparametrization.
Akaike's information criterion
formula is given by
\[
\operatorname{AIC}=n\operatorname{ln}(\operatorname
{RMSS})+2(p+1)-n\operatorname{ln}(n).
\]
Specially, we have the following three AICs for the above chosen
three models.
\begin{eqnarray*}
\operatorname{AIC}_u&=&n\operatorname{ln}\bigl(\operatorname{tr}\bigl((
Y-X\widehat
\Theta
_uZ_1')'( Y-X\widehat\Theta
_uZ_1')\bigr)\bigr)+2(p_u+1)-n\operatorname{ln}(n),
\\
\operatorname{AIC}_o&=&n\operatorname{ln}\bigl(\operatorname{tr}\bigl((
Y-X\widehat
\Theta
_oZ_2')'( Y-X\widehat\Theta
_oZ_2')\bigr)\bigr)+2(p_o+1)-n\operatorname{ln}(n)
\end{eqnarray*}
and
\[
\operatorname{AIC}_a=n\operatorname{ln}\bigl(\operatorname{tr}\bigl((
Y-X_1\widehat\Theta
_1Z_1'-X_2\widehat\Theta_2Z_2')'( Y-X_1\widehat\Theta
_1Z_1'-X_2\widehat\Theta
_2Z_2')\bigr)\bigr)+2(p_a+1)-n\operatorname{ln}(n),\vadjust{\goodbreak}
\]
where $p_u$, $p_o$ are the numbers of parameters in $\Theta_u$,
$\Theta_o$, respectively, and $p_a$ is sum of the numbers of
parameters in $\Theta_1$ and $\Theta_2$.

In our simulation, consider $n_1=n_2=n/2$, replication times
$N=10~000$ and $\rho=0.2, 0.5$ and~$0.8$, respectively.

With $N$ replication times, the average values of $\operatorname{AIC}_u$,
$\operatorname{AIC}_o$ and $\operatorname{AIC}_a$ are denoted by
\[
\operatorname{AIC}(\psi_u)= \frac1N \sum_{i=1}^N\operatorname{AIC}_u^i,\qquad
\operatorname{AIC}(\psi_o)= \frac1N \sum_{i=1}^N\operatorname
{AIC}_o^i,\qquad \operatorname
{AIC}(\psi_a)= \frac1N
\sum_{i=1}^N\operatorname{AIC}_a^i.
\]

%
%
\begin{figure}

\includegraphics{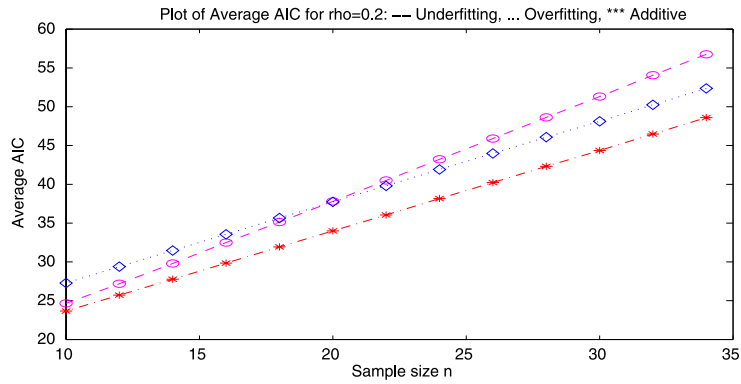}

\caption{AIC($\psi_u$), AIC($\psi_o$), AIC($\psi_a$)
and sample size n for $\rho=0.2$.}\label{fig:1}
\end{figure}

%
\begin{figure}[b]

\includegraphics{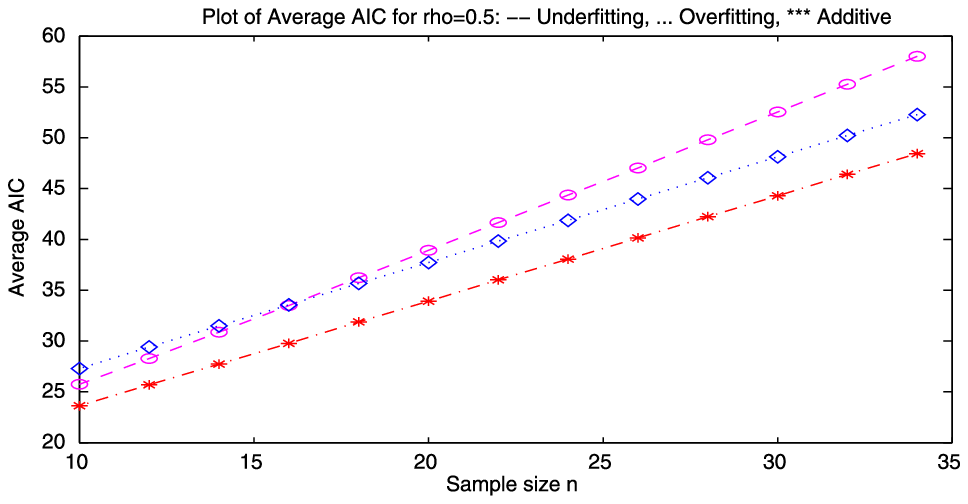}

\caption{AIC($\psi_u$), AIC($\psi_o$), AIC($\psi_a$)
and sample size n for $\rho=0.5$.}\label{fig:2}
\end{figure}

%
\begin{figure}

\includegraphics{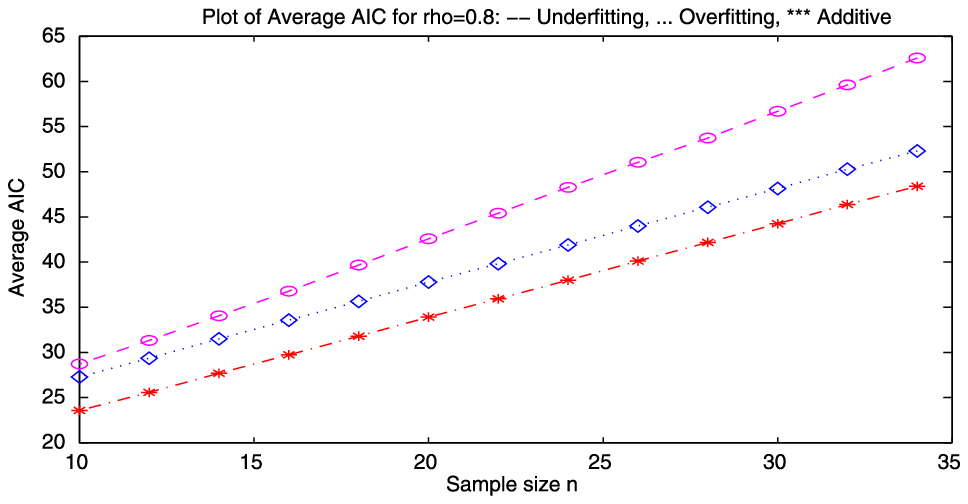}

\caption{AIC($\psi_u$), AIC($\psi_o$), AIC($\psi_a$)
and sample size n for $\rho=0.8$.}\label{fig:3}
\end{figure}

The relations between the sample size $n$ and $\operatorname{AIC}(\psi_u)$,
$\operatorname{AIC}(\psi_o)$ and $\operatorname{AIC}(\psi_a)$ are
illustrated in Figures
1--3 for $\rho=0.2,
0.5$ and $0.8$, respectively. We can make the following conclusions
from these curves:

\begin{longlist}
\item[(1)] Akaike's information criterion of the true model $\psi_a$ remains
to be uniformly smallest
for all cases of $\rho=0.2, 0.5$ and $0.8$. The trend becomes particularly
obvious as the sample size $n$ increases. We believe that the
conclusion is true for all $\rho\in(0, 1)$.

\item[(2)] The curve for AIC of the true model $\psi_a$ and the curve for
AIC of the overfitted model~$\psi_o$ are parallel. It implies that
the difference between AIC for the true model $\psi_a$ and AIC of the
overfitted model $\psi_o$ is a constant. The constant is due to the
penalty for overparametrization. This shows that it is not
significant for the difference between the RMSS for the true model
$\psi
_a$ and the RMSS for the overfitted model. Overfitting gets a penalty for
overparametrization and leads to a bigger AIC.

\item[(3)] The underfitted model $\psi_u$ seems to have a bigger AIC than
the overfitted model. It means that underfitting incurs more loss
than overfitting does in the terms of AIC. The loss becomes larger
and larger as the sample size increases or $\rho$ is closer and
closer to 1.

\item[(4)] The curve for AIC of the underfitted model becomes a little bit
steeper as
$\rho$ is gradually close to 0, while the curve for AIC of the
overfitted model and the curve for AIC of the true model seem to be
unrelated to $\rho$.
\end{longlist}

In conclusion, using the additive growth curve model (\ref{e11}) with
orthogonal design matrices has an obvious advantage over using the
traditional growth curve model in model specification and then in
parameter estimation.

\section{A numerical example}\label{sec6}

The numerical example, stated in \cite{Potthoff64}, about a
certain measurement in a dental study on 11 girls and 16 boys at 4
different ages -- 8, 10, 12 and 14 -- is employed here (see Table~\ref{tab1}) to
illustrate the ideas and techniques stated in the paper.
%
\begin{table}
\caption{Measurements on 11 girls and 16 boys,
at 4 different ages --$8$, $10$, $12$, $14$}\label{tab1}
\begin{tabular*}{\textwidth}{@{\extracolsep{\fill}}ld{2.2}d{2.2}d{2.2}d{2.2}ld{2.2}d{2.2}d{2.2}d{2.2}@{}}
\hline
\multicolumn{1}{@{}l}{Girls} & \multicolumn{1}{l}{8} & \multicolumn{1}{l}{10} &
\multicolumn{1}{l}{12} & \multicolumn{1}{l}{14} &
\multicolumn{1}{l}{Boys} &
\multicolumn{1}{l}{8}&\multicolumn{1}{l}{10}&
\multicolumn{1}{l}{12}&\multicolumn{1}{l@{}}{14} \\
\hline
\phantom{0}1& 21 & 20 & 21.5 & 23 & \phantom{0}1 & 26 & 25 & 29 & 31 \\
\phantom{0}2& 21 & 21.5 & 24 & 25.5 & \phantom{0}2 & 21.5 & 22.5 & 23 & 26.5 \\
\phantom{0}3& 20.5 & 24 & 24.5 & 26 & \phantom{0}3 & 23 & 22.5 & 24 & 27.5 \\
\phantom{0}4& 23.5 & 24.5 & 25 & 26.5 & \phantom{0}4 & 25.5 & 27.5 & 26.5 & 27 \\
\phantom{0}5& 21.5 & 23 & 22.5 & 23.5 & \phantom{0}5 & 20 & 23.5 & 22.5 & 26 \\
\phantom{0}6& 20 & 21 & 21 & 22.5 & \phantom{0}6 & 24.5 & 25.5 & 27 & 28.5 \\
\phantom{0}7& 21.5 & 22.5 & 23 & 25 & \phantom{0}7 & 22 & 22 & 24.5 & 26.5 \\
\phantom{0}8& 23 & 23 & 23.5 & 24 & \phantom{0}8 & 24 & 21.5 & 24.5 & 25.5 \\
\phantom{0}9& 20 & 21 & 22 & 21.5 & \phantom{0}9 & 23 & 20.5 & 31 & 26 \\
10& 16.5 & 19 & 19 & 19.5 & 10 & 27.5 & 28 & 31 & 31.5 \\
11& 24.5 & 25 & 28 & 28 & 11 & 23 & 23 & 23.5 & 25 \\
& & & & & 12 & 21.5 & 23.5 & 24 & 28 \\
& & & & & 13 & 17 & 24.5 & 26 & 29.5 \\
& & & & & 14 & 22.5 & 25.5 & 25.5 & 26 \\
& & & & & 15 & 23 & 24.5 & 26 & 30 \\
& & & & & 16 & 22 & 21.5 & 23.5 & 25 \\[6pt]
Mean & 21.18 & 22.23 & 23.09 & 24.09 &
Mean &
22.87 &
23.81 & 25.72 & 27.47 \\
\hline
\end{tabular*}
\end{table}

Prior to making the model specification, we do not know whether the
distances, in millimeters, from the center of the pituitary to the
pteryo-maxillary fissure of these girls and boys follow two
polynomial functions of time $t$ with a same degree. So we assume
that the distances for girls and for boys follow two polynomial
functions of time $t$ with different degrees $g$ and $b$ (set $1\leq
g, b\leq3$).

Based on the model $(\ref{e11})$, we think of these observations
as realizations of the following model:
\[
Y=X_1\Theta_1Z_g'+X_2\Theta_2Z_b'+\mathcal E,
\]
where
\[
X_1=
\pmatrix{\mathbf1_{11}\cr
\mathbf{0}},\qquad
\Theta_1=
\pmatrix{\theta_{11} & \cdots&
\theta_{1g}},\qquad Z_g'=
\pmatrix{1 & 1& 1 & 1\cr
\cdot& \cdot& \cdot& \cdot\cr
t_1^{g}& t_2^{g}& t_3^{g}&
t_{4}^{g}}\qquad
\mbox{for }1\leq g\leq3,
\]
and
\[
X_2=
\pmatrix{\mathbf{0}\cr
\mathbf1_{16}},\qquad
\Theta_2=
\pmatrix{\theta_{21} & \cdots&
\theta_{2b}}, \qquad Z_b'=
\pmatrix{1 & 1& 1 & 1\cr
\cdot& \cdot& \cdot& \cdot\cr
t_1^{b}& t_2^{b}& t_3^{b}& t_{4}^{b}}
\qquad  \mbox{for }1\leq b\leq3.
\]

We should trade the effect of the RMSS
from the simple ``true'' model and the loss from overparameterization.
Due to setting $1\leq g,b\leq3$, we can structure nine models for
selection. The corresponding AICs of the nine models are displayed in
Table \ref{tab2}.

The best model is the one with the minimum AIC. Based on AIC, the model
with parameter pairs $(1,1)$ is
best, that is, the growth curves for girls and boys are two linear
equations of time~$t$. Our conclusion of model specification is
consistent with the chosen model of~\cite{Potthoff64}.

\begin{table}
\caption{Parameter pair $(g,b)$ and AICs for 9
models}\label{tab2}
\begin{tabular*}{\textwidth}{@{\extracolsep{\fill}}llllll@{}}
\hline
($g$, $b$) & AIC & ($g$, $d$) & AIC & ($g$, $b$) & AIC\\
\hline
(1, 1) & 90.4011$^*$& (1, 2) & 92.2497 & (1, 3) & 94.1817 \\
(2, 1) & 92.4009& (2, 2) & 94.2495 & (2, 3) & 96.1815  \\
(3, 1) & 94.3972 & (3, 2) & 96.2458 & (3, 3) & 98.1777 \\
\hline
\end{tabular*}
\end{table}
%

\section{Concluding remarks}\label{sec7}

When observations of a repeated measurement at multiple time points
follow polynomial functions with different degrees rather than the
same degree, using the traditional growth curve model may cause
underfitting or overfitting. To
avoid these troubles, we proposed an additive growth curve model (\ref
{e11}) with orthogonal design matrices that allows us to fit the
data and then estimate parameters in a more parsimonious and
efficient way than the traditional growth curve model. Obviously,
the proposed additive growth curve model can not be included in the
extended growth curve models investigated by Kollo and von Rosen
\cite{Kollo05}, Chapter~4,
and Verbyla and Venables \cite{Verbyla88}.

In the paper, we explored the least-squares approach to derive
two-stage GLS estimators for the regression
coefficients, where an invariant and unbiased quadratic estimator for
the covariance of observations is taken as the first-stage estimator. We
investigated the properties of these estimators, including
unbiasedness, consistency and asymptotic normality.

Simulation studies and a numerical example are given to illustrate
the efficiency and parsimony of the proposed model for model
specification in the sense of minimizing AIC compared to the
traditional growth curve model. It
follows that our additive growth curve model and the least-squares
estimation for regression coefficients are competitive alternatives
to the traditional growth curve model.

\begin{appendix}
\section*{Appendix}\label{app}
\begin{pf*}{Proof of Theorem \protect\ref{l21}}
Put $T=(X_1\otimes
Z_1,\ldots,X_k\otimes Z_k)$ and $\bolds{\beta}=((\operatorname
{vec}(\Theta_1))',
\ldots, (\operatorname{vec}(\Theta_k))')'$. Then the model (\ref
{e11}) can
be rewritten as $\operatorname{vec}(\bolds{\mu})=T\bolds{\beta}$
where $\mathscr
C(T)=\mathscr C(X_1\otimes Z_1)+\cdots+\mathscr C(X_k\otimes Z_k)$.

(a) To obtain the two-stage generalized least square estimate of
$\bolds{\mu}$ from the data $Y$ is equivalent to applying the ordinary
least square method to the following model
\renewcommand{\theequation}{\arabic{equation}}
\setcounter{equation}{14}
\begin{equation}\label{e29}
\operatorname{vec}(Z)=\operatorname{vec}(\bolds{\nu})+\bigl(I\otimes
\Sigma^{-1/2}(Y)\bigr)\operatorname{vec}(\mathcal{E}),
\end{equation}
where $\operatorname{vec}(Z)=(I\otimes\Sigma^{-1/2}(Y))\operatorname
{vec}(Y)$ and
$\operatorname{vec}(\bolds{\nu})=(I\otimes\Sigma^{-1/2}(Y))T\bolds
{\beta}$. So the
ordinary least square estimator of $\operatorname{vec}(\bolds{\nu})$ is
\begin{equation}\label{e210}
\operatorname{vec}(\hat{\bolds{\nu}}_{\mathrm{ols}}(Z)=P_{(I\otimes\Sigma
^{-1/2}(Y))T}
\operatorname{vec}(Z).
\end{equation}
Thus by equations (\ref{e29}) and (\ref{e210})
\begin{equation}\label{e211}
\operatorname{vec}(\hat{\bolds{\mu}}(Y))=\bigl(I\otimes\Sigma
^{1/2}(Y)\bigr)P_{(I\otimes
\Sigma^{-1/2}(Y))T}\bigl(I\otimes\Sigma^{-1/2}(Y)\bigr)\operatorname{vec}(Y).
\end{equation}
Since
\begin{equation}\label{e212}
P_{(I\otimes\Sigma^{-1/2}(Y))T} =\bigl(I\otimes
\Sigma^{-1/2}(Y)\bigr)T\bigl(T' \bigl(I\otimes\Sigma(Y)\bigr)^{-1}T\bigr) ^+T' \bigl(I\otimes
\Sigma^{-1/2}(Y)\bigr),
\end{equation}
by equations (\ref{e211}) and (\ref{e212}), we obtain
\[
\operatorname{vec}(\hat{\bolds{\mu}}(Y))=T\bigl(T' \bigl(I\otimes\Sigma
(Y)\bigr)^{-1}T\bigr)^+T'
\bigl(I\otimes\Sigma(Y)\bigr)^{-1}\operatorname{vec}(Y).
\]
By Kronecker product operations and (\ref{e27}),
$\operatorname{vec}(\hat{\bolds{\mu}}(Y))$ reduces to
\begin{equation}\label{e214}
\operatorname{vec}(\hat{\bolds{\mu}}(Y))=\sum_{i=1}^k\{
X_i(X_i'X_i)^-X_i'\otimes
Z_i(Z_i'\Sigma^{-1}(Y)Z_i)^+Z_i'\Sigma^{-1}(Y)\}\operatorname{vec}(Y).
\end{equation}
Since $(Z_i(Z_i'Z_i)^{-}Z_i'\Sigma^{-1}(Y)Z_i(Z_i'
Z_i)^-Z_i')^+=Z_i(Z_i'\Sigma^{-1}(Y)Z_i)^+Z_i'$, in matrix language,
we obtain the expression (\ref{e27}) of $\hat{\bolds{\mu}}(Y)$ by
rewriting (\ref{e214}).

(b) It follows from equation (\ref{e21}) and the condition (\ref{e12}).

(c) To prove the unbiasedness of $\Theta_i$'s, by (\ref{e22}) and
(\ref{e27}), it suffices to show that $\hat{\bolds{\mu}}(Y)$ is an
unbiased estimator of $\bolds{\mu}$.

Since $\widehat\Sigma(Y)=\widehat\Sigma(\mathcal
E)=\widehat\Sigma(-\mathcal E)$, $H_i(-\mathcal E)=H_i(\mathcal E)$
and $\mathrm{E}(\mathcal E H_i(\mathcal E))=\mathbf{0}$. By (\ref{e25}),
$\hat{\bolds{\mu}}(Y)$ can be expressed as
\[
\hat{\bolds{\mu}}(Y)=\sum_{i=1}^kP_{X_i}Y H_i(Y) =\sum_{i=1}^k
X_i\Theta
_iZ_i'+\sum_{i=1}^kP_{X_i}{\mathcal E}H_i(\mathcal E).
\]
And
\[
\mathrm{E}(\hat{\bolds{\mu}}(Y))=\sum_{i=1}^kX_i\Theta
_iZ_i'+\sum_{i=1}^kP_{X_i}\mathrm{E}({\mathcal E}H_i(\mathcal
E))=\sum_{i=1}^kX_i\Theta_iZ_i'=\bolds{\mu},
\]
completing the proof.
\end{pf*}

\begin{pf*}{Proof of Lemma \protect\ref{l41}} Some subscript $i$'s are
ignored in the following statements. Write
$V=\frac{1}{\sqrt{n}}X_i'=[\mathbf{v}_1, \ldots, \mathbf{v}_n]$.
The transpose of $\mathbf{v}_j$ is an $m_i$-element row vector as
follows,
\[
\mathbf{v}_j'=\biggl(\frac{1}{\sqrt{n}}a_{j1}, \ldots, \frac{1}{\sqrt
{n}}a_{jm}\biggr),
\]
where $X_i =[a_{lj}]_{n\times m_i}$. By (\ref{e34}),
$VV'=n^{-1}X_i'X_i$ converges to a positive definite matrix $R_i$.
So the elements of $VV'=\mathbf{v}_1\mathbf{v_1}'+\cdots+\mathbf
{v}_n\mathbf{v}_n'$ are bounded. We claim that, for any
$j\in\{1,\ldots,n\}$, the $m_i$ elements of $\mathbf{v}_j$ are all
$\mathrm{O}(n^{-1/2})$.

If this is not true, we can assume without loss of generality that
one element of $\mathbf{v}_n$ is $\mathrm{O}(n^{p-1/2})$ with $p>0$. Then one
element of $\mathbf{v}_n\mathbf{v}_n'$ would be $\mathrm{O}(n^{2p-1})$. Hence,
the corresponding element in matrix $VV'=\mathbf{v}_1\mathbf
{v}_1'+\cdots+\mathbf{v}_{n}\mathbf{v}_{n}'$ would be $\mathrm{O}(n^{2p})$, which
is not bounded. This is a contradiction to condition (\ref{e34}).

Since
\begin{eqnarray*}
\bigl(\sqrt{n}\mathbf{s}_{i1}, \ldots, \sqrt{n}\mathbf
{s}_{in}\bigr)&=&\sqrt{n}(X_i'X_i)^{-1}X_i'=n(X_i'X_i)^{-1}\frac{1}{\sqrt{n}}X_i'\\
&=&n(X_i'X_i)^{-1}[\mathbf{v}_1, \ldots, \mathbf{v}_n],
\end{eqnarray*}
namely, for $j=1,\ldots, n$, $ \sqrt{n}\mathbf{s}_{ij}=n(X_i'
X_i)^{-1}\mathbf{v}_j$. Thus, for $j=1,\ldots,n$, the $m_i$ elements
of $\sqrt{n}\mathbf{s}_{ij}$ are also $\mathrm{O}(n^{-1/2})$, completing the
proof.
\end{pf*}

\begin{pf*}{Proof of Theorem \protect\ref{t41}}
Fix $i$. Let $\Gamma
_i=S_i\mathcal E\in\mathscr M_{m_i\times p}$. Then $\Gamma_i$ can
be rewritten as
\[
\Gamma_i=\sum_{j=1}^n\mathbf{s}_{ij}\mathcal E_{j}',
\]
where $\mathbf{s}_{ij}$ is the $j$th column vector of $\mathbf{S}_i$
and $\mathcal E_j'$ is the $j$th row vector of the matrix $\mathcal
E$ with $\mathcal E\sim\mathcal{G}(\mathbf{0}, \mathbf{I}_n \otimes
\Sigma)$.

Since $\{\mathcal E_j'\}_{j=1}^n$ are independent and identically
distributed, for $\mathbf{t}$ $\in\mathscr M_{m_i\times p}$, the
characteristic function $\Psi_n(\mathbf{t})$ of $\sqrt{n}\Gamma_i$
is given by
\begin{eqnarray*}
\Psi_n(\mathbf{t})&=&\mathrm{E}\bigl(\exp\bigl\{\mathrm{i}\operatorname
{tr}\bigl(\sqrt{n}\mathbf{t}'\Gamma
_i\bigr)\bigr\}\bigr)=\mathrm{E}\Biggl(\exp\Biggl\{\mathrm{i}\operatorname{tr}\Biggl(\sqrt{n}\mathbf{t}'\sum
_{j=1}^{n}\mathbf
{s}_{ij}\mathcal E_{j}'\Biggr)\Biggr\}\Biggr)\\
&=&\mathrm{E}\Biggl(\exp\Biggl\{\mathrm{i}\operatorname{tr}\Biggl(\sqrt{n}\sum_{j=1}^n\mathbf
{t}'\mathbf{s}_{ij}\mathcal
E_{j}'\Biggr)\Biggr\}\Biggr)=\prod_{j=1}^n\Phi\bigl(\sqrt{n}\mathbf{t}'\mathbf{s}_{ij}\bigr),
\end{eqnarray*}
where $\Phi(\cdot)$ is the characteristic function of $\mathcal
E_j'$.

Recall that for $u$ in the neighborhood of $0$,
\renewcommand{\theequation}{\arabic{equation}}
\setcounter{equation}{19}
\begin{equation}\label{e67}
\ln(1-u)=-u+f(u)\qquad \mbox{with } f(u)=\frac{1}{2}u^2+\mathrm{o}(u^2).
\end{equation}
Write $p(u)={f(u)}/{u}$, then from (\ref{e67}),
\begin{equation}\label{e68}
p(u)=\mathrm{o}(u) \qquad \mbox{as } u\rightarrow0.
\end{equation}
And
\begin{eqnarray}\label{e65}
\Phi(\mathbf{x})&=&1 -\frac{1}{2}\mathbf{x}'\Sigma\mathbf{x}
+g(\mathbf
{x})\qquad \mbox{for } \mathbf{x}\in\mathbb{R}^{m_i} \quad\mbox{and}
\nonumber
\\[-8pt]
\\[-8pt]
\nonumber
g(\mathbf{x})&=&\mathrm{o}(\parallel\mathbf{x}\parallel^2) \qquad\mbox{as } \mathbf
{x}\rightarrow
\mathbf{0}.
\end{eqnarray}
For $\varepsilon>0$, there exists $\delta(\varepsilon)>0$ such that
\begin{equation}\label{e66}
|g(\mathbf{x})|<\varepsilon{\parallel\mathbf{x}\parallel^2}
\qquad\mbox{as } 0<\parallel\mathbf{x}\parallel<\delta(\varepsilon).
\end{equation}
By (\ref{e67}) and (\ref{e65}),
\begin{eqnarray*}
\ln\bigl(\Phi\bigl(\sqrt{n} \mathbf{t}'\mathbf{s}_{ij}\bigr)\bigr)&=&\ln
\biggl(1-\frac
{n}{2}\mathbf{s}_{ij}'\mathbf{t}\Sigma\mathbf{t}'\mathbf
{s}_{ij}+g\bigl(\sqrt
{n}\mathbf{t}'\mathbf{s}_{ij}\bigr)\biggr)\\
&=&-\frac{1}{2}n\mathbf{s}_{ij}'\mathbf{t}\Sigma\mathbf{t}'\mathbf
{s}_{ij}+g\bigl(\sqrt{n}\mathbf{t}'\mathbf{s}_{ij}\bigr)
+f\biggl(\frac{1}{2}n\mathbf{s}_{ij}'\mathbf{t}\Sigma\mathbf{t}'\mathbf
{s}_{ij}-g\bigl(\sqrt{n}\mathbf{t}'\mathbf{s}_{ij}\bigr)\biggr).
\end{eqnarray*}
Therefore, the characteristic function of $\sqrt{n}\Gamma_n$ can be
decomposed as
\begin{equation}\label{e69}
\Psi_n(\mathbf{t})=\exp\Biggl\{ \sum_{j=1}^n
\ln\bigl(\Phi\bigl(\sqrt{n} \mathbf{t}'\mathbf{s}_{ij}\bigr)\bigr)\Biggr\}
\equiv\exp\biggl\{-\frac{1}{2}\alpha_n+\bolds{\beta}_n+\eta_n\biggr\},
\end{equation}
where
\begin{eqnarray*}
\alpha_n&=&\sum_{j=1}^nn\mathbf{s}_{ij}'\mathbf{t}\Sigma\mathbf
{t}'\mathbf
{s}_{ij}=\operatorname{tr}\Biggl(\sum_{j=1}^nn\mathbf{s}_{ij}'\mathbf
{t}\Sigma\mathbf
{t}'\mathbf{s}_{ij}\Biggr),
\\
\bolds{\beta}_n&=&\sum_{j=1}^ng\bigl(\sqrt{n}\mathbf{t}'\mathbf{s}_{ij}\bigr)
\end{eqnarray*}
and
\[
\eta_n= \sum_{j=1}^nf\biggl(\frac{1}{2}n\mathbf{s}_{ij}'\mathbf{t}\Sigma
\mathbf{t}'\mathbf{s}_{ij}-g\bigl(\sqrt{n}\mathbf{t}'\mathbf{s}_{ij}\bigr)\biggr).
\]
Note that
\begin{equation}\label{e699}
\sum_{j=1}^n\mathbf{s}_{ij}\mathbf{s}_{ij}'=(X_i'X_i)^{-1}.
\end{equation}
For $\alpha_n$, by (\ref{e699}), we have
\begin{equation}\label{e610}
\alpha_n=\operatorname{tr}\Biggl(n\mathbf{t}\Sigma\mathbf{t}'\sum
_{j=1}^n\mathbf
{s}_{ij}\mathbf{s}_{ij}'\Biggr)=\operatorname{tr}(n\Sigma\mathbf{t}'(X_i'X_i)^{-1}\mathbf{t}).
\end{equation}
By (\ref{e34}),
\begin{equation}\label{e611}
\lim_{n\rightarrow\infty} \alpha_n=\operatorname{vec}(\mathbf
{t})'(R_i^{-1}\otimes
\Sigma)\operatorname{vec}(\mathbf{t}).
\end{equation}
For $\bolds{\beta}_n$, by Lemma \ref{l41} and the continuity of
$\mathbf
{t}'\mathbf{s}_{ij}$, for the $\delta(\varepsilon)>0$ in (\ref{e66}),
there is an integer $N(\varepsilon)>0$ such that for
$n>N(\varepsilon)$,
\begin{equation}\label{e612}
0<\bigl\Vert\sqrt{n}\mathbf{t}'\mathbf{s}_{ij}\bigr\Vert<\delta
(\varepsilon)\qquad  \mbox{for all } j=1,\ldots,n.
\end{equation}
Take $n>N(\varepsilon)$, then by (\ref{e66}) and (\ref{e612}),
\begin{equation}\label{e613}
\bigl|g\bigl(\sqrt{n}\mathbf{t}'\mathbf{s}_{ij}\bigr)\bigr|<\bigl\Vert\sqrt{n}\mathbf
{t}'\mathbf{s}_{ij}\bigr\Vert^2\varepsilon.
\end{equation}
By (\ref{e699}),
\begin{equation}\label{e6133}
|\bolds{\beta}_n|<\sum_{j=1}^n \bigl\Vert\sqrt{n}\mathbf
{t}'\mathbf
{s}_{ij}\bigr\Vert^2\varepsilon=\varepsilon n\sum
_{j=1}^n\operatorname{tr}(\mathbf
{t}'\mathbf{s}_{ij}\mathbf{s}_{ij}'\mathbf{t})
=\varepsilon\operatorname{tr}(\mathbf{t}'n(X_iX_i)^{-1}\mathbf{t}).
\end{equation}
So by (\ref{e34}), $\limsup_{n\rightarrow\infty}|\bolds{\beta
}_n|\leq
\varepsilon\operatorname{tr}(\mathbf{t}'R_i^{-1} \mathbf{t})$.
Since $\varepsilon>0
$ is arbitrary, we obtain
\begin{equation}\label{e614}
\lim_{n \rightarrow\infty}\bolds{\beta}_n=0.
\end{equation}
And for $\eta_n$, let
\[
\lambda_j= \frac{1}{2}\bigl(\sqrt{n}\mathbf{t}'\mathbf{s}_{ij}\bigr)'\Sigma
\bigl(\sqrt
{n}\mathbf{t}'\mathbf{s}_{ij}\bigr)-g\bigl(\sqrt{n}\mathbf{t}'\mathbf{s}_{ij}\bigr).
\]
Thus, by (\ref{e613}),
\begin{equation}\label{e615}
|\lambda_{j}|<\frac{1}{2}\bigl(\sqrt{n}\mathbf{t}'\mathbf
{s}_{ij}\bigr)'\Sigma
\bigl(\sqrt{n}\mathbf{t}'\mathbf{s}_{ij}\bigr)+\bigl\Vert\sqrt{n}\mathbf
{t}'\mathbf
{s}_{ij}\bigr\Vert^2\varepsilon.
\end{equation}
Take $n>N(\varepsilon)$, by Lemma \ref{l41}, the continuity of
$\mathbf{t}'\mathbf{s}_{ij}$ and (\ref{e68}), increasing
$N(\varepsilon)$ if necessary, we may suppose that for all $j$,
$|p(\lambda_{j})|<\varepsilon$. Since
$f(\lambda_{j})=p(\lambda_{j})\lambda_{j}$,
\[
|\eta_n|\leq\sum_{j=1}^n|f(\lambda_j)|=\sum_{j=1}^n|p(\lambda
_j)||\lambda_j| \leq
\sum_{j=1}^n\varepsilon|\lambda_{j}|.
\]
By (\ref{e615}),
\[
|\eta_n|\leq\sum_{j=1}^n\biggl[\frac{\varepsilon}{2}
\operatorname{tr}\bigl(\sqrt{n}\mathbf{s}_{ij}'\mathbf{t}\Sigma\mathbf
{t}'\sqrt{n}\mathbf
{s}_{ij}\bigr)+\bigl\Vert\sqrt{n}\mathbf{t}\mathbf{s}_{ij}\bigr\Vert
^2\varepsilon^2\biggr].
\]
Then, taking the same operations as (\ref{e610}) and (\ref
{e6133}), we obtain the following inequality
\[
|\eta_n|\leq\biggl[\frac{\varepsilon}{2}\operatorname{tr}
(n(X_i'X_i)^{-1}\mathbf{t}\Sigma\mathbf{t}')
+\varepsilon^2\operatorname{tr}(\mathbf{t}'n(X_i'X_i)^{-1}\mathbf{t})\biggr],
\]
namely, by (\ref{e34})
\begin{equation}\label{e616}
|\eta_n|\leq\frac{\varepsilon}{2}\operatorname{tr}\Biggl(\sum_{l=1}^i
R_i^{-1}\mathbf{t}\Sigma\mathbf{t}'\Biggr)
+\varepsilon^2\operatorname{tr}(\mathbf{t}'R_i^{-1}\mathbf{t}).
\end{equation}
Due to arbitrary of $\varepsilon$ and (\ref{e616}),
\begin{equation}\label{e617}
\lim_{n \rightarrow\infty}\eta_n=0.
\end{equation}
By (\ref{e611}), (\ref{e614}) and (\ref{e617}), we obtain from
(\ref{e69}),
\begin{equation}\label{e618}
\lim_{n \rightarrow\infty}\Psi_n(\mathbf{t})=\exp\biggl\{-\frac
{1}{2}\operatorname{vec}
(\mathbf{t})'(R_i^{-1}\otimes
\Sigma)\operatorname{vec}(\mathbf{t}) \biggr\}.
\end{equation}
So by Levy's continuity theorem, $\sqrt{n}\Gamma_i $ converges in
distribution to $\mathcal N_{m_i\times p}(\mathbf{0}, R_i^{-1}\otimes
\Sigma)$, completing the proof of the desired result.
\end{pf*}
\end{appendix}

\section*{Acknowledgements}

Hu's research was supported by the National Natural Science Foundation
of China (NSFC) Grants 10971126. Shanghai University of Finance and
Economics provides partial funding for
Hu's work and You's work through Project 211, Phase III, and Shanghai
Leading Academic Discipline Project B803. Yan's research was funded by
the Natural Sciences and Engineering Research Council of Canada, RGPIN
371505-10.
The authors are deeply grateful to the anonymous referee and the
associate editor for their helpful comments, which led to the
greatly improved version of this paper.

\printhistory

\end{document}